\documentclass[11pt]{article}
\usepackage[utf8]{inputenc}

\usepackage{geometry}                
\usepackage{graphicx}
\usepackage{amssymb}
\usepackage{amsthm}
\usepackage{amsmath}

\usepackage[title]{appendix}%

\usepackage{natbib}
\usepackage[font=small]{caption}
\usepackage{dcolumn}
\usepackage{multirow,booktabs}

\usepackage{float}
\newcolumntype{L}[1]{>{\arraybackslash}p{#1}}
\newcolumntype{K}[1]{>{\centering\arraybackslash}p{#1}}
\newcolumntype{R}[1]{>{\raggedleft\arraybackslash}p{#1}}
\usepackage{longtable}

\newtheorem{lemma}{Lemma}
\newtheorem{prop}{Proposition}

\usepackage{epstopdf}
\begin{document}

\title{Multi-period Newsvendor Model}
\author{Valentyn Khokhlov}
\date{February 10, 2026}

\maketitle

\begin{abstract}
The newsvendor model is a well-known stochastic model for inventory management; however, it was originally developed for a single-period context and focuses on trading companies. This paper proposes an extension of the newsvendor model into a mutli-period setting, aiming to develop a decision-making tool for manufacturing firms to determine the optimal production batch size. The objective function is to maximize operating profit in accordance with generally accepted accounting principles. The model can also incorporate overhead costs, such as warehousing, shrinkage, cost of capital, and lead time between the production decision and output. Monte Carlo simulations demonstrate that the proposed model results in higher profitability compared to other newsvendor models used in our analysis, as well as the safety stock buffer approach. The key feature explaining its outperformance is better adaptability of the production batch size, that leads to fewer stock-outs relative to other newsvendor models and lower inventory levels compared to the safety stock buffer approach. The robustness analysis shows that the proposed model is quite tolerant of mismatches between the ``model" and the ``true" demand distributions. Finally, we provide some recommendations on selecting  the appropriate ``model" distribution for different SKUs.

\paragraph{Keywords} inventory management, newsvendor model, newsboy model, optimal order size, optimal production batch size   
\end{abstract}

\section{Introduction}

The newsvendor model (also known as the newsboy model) is a well-known single-period stochastic model for inventory management. An extensive literary review of this model was presented in \cite{Khouja1999}. Although the model was developed primarily for a trading company and addresses the problem of determining the optimal inventory level to balance the profits lost due to stock shortages against the costs of holding excess inventory, it can also be extended to cover the case of a manufacturing company seeking optimal production quantities under uncertain demand. A significant body of research has been dedicated to the newsvendor model, and we will provide a brief review of the major developments in the following section.  However, \cite{Ward1991} has noted that this theoretical model often does not meet practical needs of decision-making. In our opinion, one of the key shortcomings of the model is its single-period nature, which is hardly a plausible assumption for a running business.

The original newsvendor problem models the decision-making of a newsboy who earns money by selling today's newspapers. He must determine the quantity of newspapers to purchase early in the morning and then sell throughout the day facing uncertain demand. If demand exceeds his stock, he incurs the lost margin that he could have earned. Conversely, if at the end of the day he is left with some unsold newspapers, that stock becomes useless and the money paid for it is lost. It is the latter assumption that justifies the single-period nature of the newsvendor model. In most other businesses, unsold inventory does not become immediately obsolete. However, as we show in our literary review, most of  researchers have not challenged this assumption of the original model. 

The idea for this paper emerged when addressing a real-world inventory management problem at a clothing manufacturer that currently uses a safety-buffer approach with an optimal reorder point. The newsvendor model seemed a good starting point, but its single-period nature immediately rendered it unusable, so the need for a mutli-period model arose. Additionally, we observed another problem with state-of-the-art newvendor model developments having become increasingly cumbersome and mathematically sophisticated. Factory managers are not the persons who have taken advanced calculus courses, so they require models that are understandable to them and their ERP developers, ideally rendering closed-form optimal solutions. Software engineers focus their efforts on information exchange, business process automation, user-friendliness, cybersecurity, but not on solving systems of differential equations or doing complex mathematical optimization. Therefore, the model should avoid such complexity if we want it to be practically useful. Factory managers tend to value robustness over sophistication and flexibility; for example, a model that is robust to distributional mismatches or moderate input estimation errors is much more likely to be implemented than one that has multiple factors and many inputs, but is highly sensitive to their specification or estimation.

In this paper, we propose an approach for developing a multi-period newsvendor model that meets the above-mentioned criteria, test its effectiveness and robustness using Monte-Carlo simulation based on actual production and demand data from the clothing manufacturer, compare and contrast the proposed multi-period model with the original newsvendor models and with the current safety-buffer approach used by the manufacturer. The paper is structured as follows: Section 2 offers a brief literature review on the newsvendor model; Section 3 presents the formulation of the original model and its extension for a manufacturing company; in Section 4 we develop a multi-period model and discuss its practical implications; Section 5 describes the simulation methodology and presents numerical results; in Section 6 we investigate the robustness of the models  when the assumed, or ``model", demand distribution differs from the actual, or ``true", one; finally, we offer conclusions and suggestions in Section 7.

\section{Literature Review}

The problem of determining optimal inventory levels or production quantities has long been a focus of researchers. One of the first formulations of the newsvendor problem was proposed by \cite{Arrow1951}. The model was extended by \cite{Kabak1978} to include secondary vendors. A major development was introduced by \cite{Ismail1979}  with their cost-volume-profitability (CVP) model, which was refined and extended in \cite{Lau1980}, \cite{Norland1980} and so on. 

A popular research topic has been adding backordering capabilities to the newsvendor model to address shortages during the sale period. It was pioneered by \cite{Gallego1993} and \cite{Khouja1996}, and subsequent developments, such as \cite{Lodree2007}, \cite{Lodree2008}, \cite{Lee2010}, have incorporated lost sales, lost contracts, customer impatience into the model. Some recent works following this direction include \cite{Lee2020}, \cite{Zhang2023}, \cite{Pando2024_amm}. 

Another popular area for extending the newsvenor model is moving from a fixed-price assumption to a price-dependent demand. While the original model has one decision variable, the amount of goods to purchase or produce, extended models add a second one, the sales price. \cite{Petruzzi1999} first incorporated a price-dependent demand function into the newsvendor model, then \cite{Ahmadi2015} added dynamic pricing capabilities with multiple demand classes, \cite{Singer2021} included price-sensitive demand in stochastic production problems, \cite{Hrabec2023} researched the effects of price on demand variance.

Our review of existing research shows that relatively few researchers have investigated the newsvendor model in multi-period context. One notable paper by \cite{Ullah2019} uses a distribution-free approach to derive the optimal order quantity while accounting for both leftover inventory and price-dependent stochastic demand. However, there are clear differences from our research. \cite{Ullah2019} primarily focus on trading companies, as opposed to manufacturing firms in our research; their model does not adhere to generally accepted accounting principles, and their distribution-free approach has shortcomings when dealing with goods with specific demand distributions. Therefore, our literature review highlights gaps in current research on the newsvendor model and justifies the need for further development in a multi-period context.

\section{The Original Newsvendor Model}

\textbf{The original model (Model 1).} This model assumes that the goods are separable, planning is performed over a single period, and the order is placed such that the goods are available on the initial date. Let the demand be represented by a random variable $D$ that follows a probability distribution with cumulative distribution function $F(x)$ and probability density function $f(x)$. The cost per unit of goods sold is $c$, and the selling price is fixed at $p$. All unsold goods at the end of the period are forfeited. Under these assumptions, the expected profit the newsboy earns during the period is
\begin{equation}E[\pi] = pE[\min\{D,q\}] - cq \label{profit_1}\end{equation}
where $q$ is the order quantity (the decision variable). 

Note that $E[\min\{D,q\}]$ represents the expected amount of the goods sold. The optimal solution is based on the following lemma, which we include here, as the proposed model will also utilize its proof.

\begin{lemma}
\label{lemma}
$$ \frac{\partial}{\partial q} E[\min\{D,q\}] = 1 - F(q)  $$
\begin{proof}
We can expand the expected amount of the goods sold as 
\begin{align*}  
E[\min\{D,q\}] &= E[D|D<q]P[D<q] + E[q|D\geq q]P[D\geq q]  \\
                                 &= \frac{\int^{q}_{-\infty}{xf(x)dx}}{F(q)}F(q) + q(1-F(q)) \\
                                 &= \int^{q}_{-\infty}{xf(x)dx} + q-qF(q)
\end{align*}
and then take its partial derivative with respect to $q$:
\begin{align*}  
\frac{\partial}{\partial q} E[\min\{D,q\}] &= qf(q) + 1 - F(q) - qf(q) = 1 - F(q)
\end{align*}
\end{proof}
\end{lemma}

The optimal solution to \eqref{profit_1} follows immediately from this lemma:
\begin{equation} q^* = F^{-1}\Big(\frac{p-c}{p}\Big) \label{opt_1} \end{equation}

\textbf{The extended model (Model 2).} Although the original model targets a reseller, its extension has been developed to account for the specifics of a manufacturing company that has fixed costs $c_f$, variable costs $c_v$ (per unit produced), inventory holding costs $h$ (per unit of inventory), and all other variables remain the same as in Model 1. Under the same assumptions as in the original model, the expected profit can be expressed as
\begin{equation}
E[\pi] = pE[\min\{D,Q\}] - c_f - c_v(Q-q_0)  - hE[\max\{Q-D,0\}] 
\label{profit_2}
\end{equation}
where $Q$ is the total inventory quantity (the decision variable) and $q_0$ is the inventory available at the beginning of the period. 

The holding costs are accrued to the amount of the leftover inventory at the end of the period, which can be expressed as $E[\max\{Q-D,0\}]=Q-E[\min\{D,Q\}]$. Using Lemma \ref{lemma}, we can derive the optimal solution for \eqref{profit_2}:
\begin{equation}
    Q^* = F^{-1}\Big(\frac{p-c_v}{p+h}\Big)
    \label{opt_2_orig}
\end{equation}
In most practical cases $h << p$, so it may appear that this solution does not materially differ from \eqref{opt_1}. However, unlike the original newsvendor model, this case requires subtracting the initial inventory, which adds multi-period capabilities to the model. The optimal amount of goods to be produced (ordered) for the current period $q^*$ can therefore be calculated as
\begin{equation}
    q^* = Q^* - q_0 = F^{-1}\Big(\frac{p-c_v}{p+h}\Big) - q_0
\label{opt_2}    
\end{equation}
A more careful examination of Model 2 shows that it does not adhere to the generally accepted accounting principle of matching costs with revenues. Specifically, the model recognizes revenues for the $E[\min\{D,Q\}]$ items sold during the period, but variable costs are accrued for the $q=Q-q_0$ items produced. So this model, as well as any models that share this shortcoming (e.g. the model developed by \cite{Ullah2019}), are not compatible with current accounting standards. 

Model 2 with its optimal production quantity \eqref{opt_2} may appear to be a working multi-period model even though it does not follow accounting principles. However, as our numeric simulations demonstrate (details on the methodology are provided in Section 5), this model simply does not perform as expected. When calculating the optimal production quantity $q^*$ with \eqref{opt_2}, we have always ended up with inferior profitability compared to \eqref{opt_2_orig} or \eqref{opt_1} (the two being practically indistinguishable when $h << p$) primarily due to excessive stock-outs.

\section{The Multi-Period Newsvendor Model}

\textbf{The improved model (Model 3).} The multi-period newsvendor model proposed in this paper represents an improvement over Model 2 and is based on the following assumptions:
\begin{itemize}
    \item The planning is performed for the next production period within a multi-period environment, where leftover inventory can be carried over from one period to the next
    \item Only one batch can be produced for each period, with production completed just before the start of the period
    \item The selling price and the costs are constant within the planning period, although they can vary from period to period (i.e. for seasonality, sales and so on)
    \item There is a lead time for production (more details provided in the \textbf{Lead time considerations} subsection below)
    \item Inventory holding costs---incorporating warehousing costs, shrinkage costs, and the cost of capital tied up in inventories---are proportional to the amount of inventory
\end{itemize}

The final assumption may require clarification. While there are various approaches for estimating the components of inventory holding costs, in most cases these costs are allocated either to the quantity of the goods, the space they occupy, or their cost. All these measures are proportional to the amount of inventory; therefore, we do not need  model specific factors separately, but combine them into a single constant $h$. 

Under these assumptions, the expected profits over the next planning period will be
\begin{equation}
E[\pi] = (p-c_v)E[\min\{D,Q\}] - c_f - h\frac{Q+E[\max\{Q-D,0\}]}{2}
\label{profit_3}
\end{equation}
where $Q$ is the amount of inventory that should be available at the beginning of the planning period (a sum of the inventory carried over from the previous period, $\hat{q}_0$, and the produces quantity, $q$), $D$ is the demand (a random variable), $p$ is the selling price during the planning period, $c_v$ is the cost of goods sold (see the \textbf{Costing considerations} subsection below), $c_f$ is the fixed periodic cost, $h$ is the inventory holding cost.

The objective function \eqref{profit_3} is quite similar to \eqref{profit_2}, with only one difference: we match variable costs to revenues rather than accruing them to the current period's production. Another minor change is that we calculate inventory holding costs not for the final inventory quantity, but for the average amount of inventory held during the period. Assuming that inventory is depleted at a constant rate from the initial quantity $Q$ to the expected final amount $E[\max\{Q-D,0\}]$, the average inventory held during the planning period is $\frac{Q+E[\max\{Q-D,0\}]}{2}$.

\textbf{Lead time considerations.} In Model 2 the term $q_0$ represented the inventory available at the beginning of the period, at the same time its value was used to calculate $q^*$, the optimal amount of goods to be produced, so $q_0$ has to be known before the current-period order is placed. Thus, an implicit assumption was made that the new order would be available immediately, allowing us both to know $q_0$ and to calculate $q^*$ on the same date. In practice, however, there is some lead time between placing an order and receiving the goods. We must make a production decision and place an order at least $T_l$ days before the start of the current period, where $T_l$ is the lead time duration in days. Therefore, the amount of inventory carrier over from the previous period is uncertain at the date the order size decision is made, which is why in Model 3 we denote its estimated value as $\hat{q}_0$ instead of using a certain value of $q_0$. 

A simple and practical way to estimate $\hat{q}_0$ is to use the last known inventory level at the date when the order size decision is made (we will keep denoting it as $q_0$) and subtract the expected demand until the end of the period. If we place the order exactly $T_l$ days before the end of the period and the period length is $T$, then
\begin{equation}
    \hat{q}_0 = q_0 - E[D]\frac{T_l}{T}
    \label{q_hat}
\end{equation}

\textbf{Costing considerations.} The original models do not specifically focus on costing details. Model 1 treats all costs as variable, while Model 2 splits them into variable and fixed costs. The latter notation is used in our model, but we can refine or adjust the meanings of the variables $c_v$ and $c_f$ in \eqref{profit_3}. Typically, managerial accounting distinguishes between variable and fixed costs, but for our model $c_v$ should only include those variable costs that are accrued to inventory (and transferred from period to period), while all other costs, whether fixed or variable, should be incorporated into $c_f$. 

This approach aligns more closely with the generally accepted accounting principles, as we can directly allocate the costs of goods sold per unit to $c_v$, and the periodic costs to $c_f$. This particular allocation is most suitable for companies that use  average inventory costing, as $c_v$ is this case would be the average cost of goods produced—--both historically and within the current period. For companies that employ different inventory costing methods, FIFO or LIFO, there would be no single $c_v$ and  different costs would have to be assigned to different inventory layers. For example, under FIFO, the expression $(p-c_v)E[\min\{D,Q\}]$ becomes much more cumbersome:
$$ pE[\min\{D,Q\}] - \sum_i{c_{v,i}q_i} - c_vE[\min\{D,Q\}-\sum_i{q_i}]$$
where $i$ indexes the past inventory layers, $c_{v,i}$ is the cost of goods in layer $i$, $q_i$ is the amount of goods in layer $i$, and we implicitly assume that the current period's production is greater than zero (because if $\sum_i{q_i} > E[D]$ we would just sell the past inventory and have not plan for additional production).

We will not pursue further developments of a FIFO- or LIFO-compatible model for two reasons. First, over a long time horizon, all methods should converge. Second, when any costing method is applied consistently, the cost of goods typically constitutes a relatively stable percentage of revenue, so we can safely assume $c_v$ to be constant in most practical scenarios. If, for some particular reason, implementing a FIFO or LIFO method becomes necessary, there is a way to incorporate it without cumbersome formulas. First, note that all past inventory layers can be combined into an estimated inventory at the end of the previous period: $\hat{q}_0 = \sum_i{q_i}$. If $\hat{q}_0 > E[D]$ then we can satisfy the next period expected demand with existing stock and $q^* = 0$. In the opposite case, we can express the effective average cost of the past inventory as $c_{v,0} = \frac{1}{\hat{q}_0}\sum_i{c_{v_i}q_i}$ and then rewrite \eqref{profit_3} in a FIFO-compatible way:
$$
E[\pi] = pE[\min\{D,Q\}] - c_{v,0}\hat{q}_0 - c_vE[\min\{D,Q\}-\hat{q}_0] - c_f - h\frac{Q+E[\max\{Q-D,0\}]}{2}
$$
The optimal solution can then be derived using the same approach as we propose below.

\begin{prop}
The optimal solution to the multi-period newsvendor model (Model 3) is
\begin{equation}
    q^* = F^{-1}\Big(\frac{p-c_v-h/2}{p-c_v+h/2}\Big)-\Big(q_0-\frac{E[D]}{T}T_l\Big)
    \label{opt_3}
\end{equation}
\begin{proof}
Considering the fact that $E[\max\{Q-D,0\}]=Q-E[\min\{D,Q\}]$ we can rewrite the objective function \eqref{profit_3} as
\begin{align*}
    E[\pi] &= (p-c_v)E[\min\{D,Q\}] - c_f - h\frac{2Q-E[\min\{D,Q\}]}{2} \\
    &= (p-c_v+h/2)E[\min\{D,Q\}] - c_f - hQ 
\end{align*}
The partial derivative of this expression with respect to $Q$, considering the results of Lemma \ref{lemma}, is
\begin{align*}
    \frac{\partial}{\partial Q}E[\pi]  &= (p-c_v+h/2)(1-F^{-1}(Q))- h
\end{align*}
so the optimal value of $Q$ would be
\begin{gather*}   
    (p-c_v+h/2)(1-F^{-1}(Q^*))- h = 0 \\
    Q^* = F^{-1}\Big(1-\frac{h}{p-c_v+h/2} \Big)=F^{-1}\Big(\frac{p-c_v-h/2}{p-c_v+h/2} \Big) 
\end{gather*} 
Considering that $q = Q - \hat{q}_0$ and using \eqref{q_hat}, we obtain \eqref{opt_3}.
\end{proof}
\end{prop}

\section{Numerical Simulation}

The numerical simulation used to evaluate the outcomes of the existing newsvendor models and the model developed in Section 4 is based on a real-world product and inventory management problem that the author encountered while consulting a clothing manufacturer. The company produced a wide range of SKUs (stock keeping units), which could be classified into two or three categories using XYZ analysis (e.g. stable/variable/sporadic demand or frequent/occasional demand). Demand was generated by  orders from the company's distributors, and the orders arrived on a daily basis during the working days (up to 23 business days per month). Production is carried out in batches with a lead time of 7 working days. Thus, we use a periodic model with a 7-day cycle.

For this simulation we selected two SKUs from different categories---SKU A with frequent demand and SKU B with occasional demand. Data for both SKUs are summarized in Table \ref{tab:SKU_params}. To simplify, we scaled both SKUs' prices to 100 and variable costs to 60. Fixed costs are allocated to SKUs on a monthly basis. The holding costs of 2.8 per period account for inventory storage, spoilage, and the cost of capital tied up in inventory on hand. The manufacturer provided daily demand data, and for SKU B it is essential that roughly 40\% of days had zero demand, so the distribution is skewed to the right. At the time, the company was using a simple inventory management approach based on a safety stock buffer, where both the reordering point and the production batch size were equal to the buffer size provided in the table.

\begin{table}[h!]
\caption{Simulation inputs for the considered SKUs}\label{tab:SKU_params}
\renewcommand{\arraystretch}{1.5}
\begin{tabular}{c c c c c K{0.5cm} K{0.6cm} c c c}
\hline 
\bfseries \multirow{2}{0.8cm}{SKU} & 
\bfseries \multirow{2}{1cm}{\centering Price, $p$} & 
\bfseries \multirow{2}{1.3cm}{\centering Variable costs,\\$c_v$} & 
\bfseries \multirow{2}{1cm}{\centering Fixed costs, $c_f$} & 
\bfseries \multirow{2}{1.3cm}{\centering Holding costs,\\$h$} & 
\multicolumn{4}{c}{\bfseries Daily Demand} &
\bfseries \multirow{2}{1cm}{\centering Safety stock buffer} \\ 
 & & & & & \bfseries min & \bfseries max & \bfseries mean & \bfseries stdev & \\ 
\hline 
A & 100 & 60 & 240,000 & 2.8 & 235 & 810 & 548.5217 & 159.3643  & 5670 \\
B & 100 & 60 &  14,000 & 2.8 &   0 &  85 &  29      &  31.28898 &  595 \\
\hline
\end{tabular}
\end{table}

Before describing the simulation methodology, we must address a particular implementation challenge. The demand data were based on daily orders (so we knew the daily demand distribution), whereas the random variable $D$ used in the newsvendor models represents periodic demand (in our case, the demand for 7 working days, which is the period length for our models). There are four possible ways to address the  discrepancy between the known daily demand distribution and the required periodic demand distribution for the models:

1. Assume the distribution of the periodic demand based on the daily demand samples and extrapolate the daily demand statistics to estimate the parameters of the periodic demand distribution. In our case, this means multiplying the expected demand by 7 and its standard deviation by $\sqrt{7}$. This approach preserves the volatility of the underlying daily data; however, it is not theoretically sound, as it implicitly assumes that the distribution of the sum of independent random variables is the same as the distribution of each individual variable.

2. Use the normal distribution, since in this particular case the sum of independent normally distributed random variables is also normally distributed. However, the demand  distributions for actual SKUs are non-normal, as they have a lower bound at zero and are typically skewed to the right. This is especially true for SKU B, for which the normal distribution is clearly inadequate.

3. Assume the distribution of the periodic demand based on the daily demand samples, but estimate its parameters using the 7-day samples. While this approach has the same theoretical flaw as the first, it allows for better adjustment of the periodic distribution parameters to fit sample statistics for the 7-day demand. While the expected value remains unaffected, the volatility is significantly lower, as this approach smooths out the daily demand fluctuations. 

4. Do not assume that the periodic demand distribution is the same as the daily demand distribution, and do not attempt to model it with a specific theoretical distribution. Thus, the only one way to calculate values of the inverse cumulative distribution function $F^{-1}(x)$ would be through a Monte Carlo simulation that generates sums of 7 random variables adhering to the daily demand distribution, with sample quantiles used to estimate the $F^{-1}$ values we require. We ended up adopting this approach.

\begin{table}[b!]
\caption{Parameters of daily demand distributions for the considered SKUs}\label{tab:distro_params}
\renewcommand{\arraystretch}{1.5}
\begin{tabular}{c | K{1cm}K{1cm} | K{1cm}K{1cm}K{2cm} | K{2cm}K{2cm}}
\hline 
\bfseries \multirow{2}{1.5cm}{\centering SKU} & 
\multicolumn{2}{c}{\bfseries Uniform} &
\multicolumn{3}{c}{\bfseries Triangular} &
\multicolumn{2}{c}{\bfseries Log-normal} \\
 & \bfseries $a$ & \bfseries $b$ & \bfseries $a$ & \bfseries $b$ & \bfseries $c$ & \bfseries $\mu_l$ & \bfseries $\sigma_l$ \\ 
\hline 
A & 235 & 810 & 235 & 810 & 600.5652 & 6.266708826 & 0.284668531 \\
B &   0 &  85 &   0 &  85 &        2 & 2.98129577  & 0.878635374 \\
\hline
\end{tabular}
\end{table}

The distributions of daily demand used in our simulation are described below, with their parameters provided in Table \ref{tab:distro_params}:

\begin{itemize}
    \setlength{\itemsep}{6pt}
    \setlength{\parskip}{0pt}
    \setlength{\parsep}{0pt}
    \item The uniform distribution. It best fits SKUs with relatively stable demand within some range $[a,b]$, where there is roughly an equal chance of observing any demand level between $a$ units and $b$ units.
    \item The triangular distribution. While not theoretically sound, this distribution is often used in practical applications when demand for a particular SKU falls within some range $[a,b]$, and the probability of observing a particular demand level grows linearly from zero at $a$ units to a maximum at $c$ units (the ``most likely" demand level), and then decreases linearly back to zero at $b$ units.
    \item The log-normal distribution. This is a theoretically sound alternative to the triangular distribution, with a slightly different shape of the probability density function. However, calculating its quantiles is much more complex in practice comparing to the previous ones. This distribution has a lower bound at zero and no upper bound. Its parameters $\mu_l, \sigma_l$ can be estimated from the sample mean $m$ and sample standard deviation $s$ as $\mu_l = ln\frac{m^2}{ \sqrt{m^2+s^2} }, \sigma_l = \sqrt{ln\big(1+\frac{s^2}{m^2}\big)}$.
\end{itemize}

Our numerical simulation is based on the Monte Carlo methodology. We perform 900 runs of simulated company operations over a period of 120 months (10 years), with each month consisting of 23 working days. The simulation is conducted separately for each of the three distributions mentioned above. During each run, we calculate the results for four different models: the safety stock buffer approach, Model 1---the original newsvendor model (\ref{opt_1}), Model 2---the extended newsvendor model (\ref{opt_2_orig}), and Model 3---the improved newsvendor model (\ref{opt_3}). The arguments of the $F^{-1}(\cdot)$ function and the optimal order sizes for each model are summarized in Table \ref{tab:arguments}. Each run proceeds as follows:

\begin{table}[b!]
\caption{The $F^{-1}(\cdot)$ argument and the optimal order size for different models}\label{tab:arguments}
\renewcommand{\arraystretch}{1.5}
\begin{tabular}{ c K{2.4cm} K{2cm} K{2.4cm} K{2.4cm} }
\hline 
\bfseries & \bfseries Safety stock & \bfseries Model 1 & \bfseries Model 2 & \bfseries Model 3 \\ 
\hline 
\bfseries $F^{-1}(\cdot)$ argument & - &  0.4 & 0.389105058 & 0.93236715 \\
\hline
& \multicolumn{4}{c}{SKU A} \\
\hline 
\bfseries Uniform    &  5670 & 3540 & 3530 & 4310$-\hat{q}_0$ \\
\bfseries Triangular &  5670 & 3705 & 3695 & 4205$-\hat{q}_0$ \\
\bfseries Log-normal &  5670 & 3715 & 3704 & 4510$-\hat{q}_0$ \\
\hline 
& \multicolumn{4}{c}{SKU B} \\
\hline 
\bfseries Uniform    &   595 &  280 &  278 &  395$-\hat{q}_0$ \\
\bfseries Triangular &   595 &  187 &  186 &  285$-\hat{q}_0$ \\
\bfseries Log-normal &   595 &  171 &  169 &  330$-\hat{q}_0$ \\
\hline
\end{tabular}
\end{table}

\begin{enumerate}
    \setlength{\itemsep}{6pt}
    \setlength{\parskip}{0pt}
    \setlength{\parsep}{0pt}
    \item We start with the inventory level equal to the optimal order size.
    \item If the current day is the day when the previous order is due, we increase the inventory level by the amount previously ordered.  
    \item For the safety stock buffer approach only: if the current inventory is below the reorder point, we place an order for the optimal order size. This effectively means we are placing orders of the same size at non-periodic points of time.
    \item For the newsvendor models: if 7 days have passed since the last order, we place a new order of the optimal order size. Since the lead time is 7 days, this
    effectively results in placing orders at periodic points of time, but their size may vary depending on the model.
    \item We simulate the demand for the current day using the selected distribution.
    \item We sell the amount of goods equal to the minimum of the daily demand and the available inventory, add the difference between the revenue and the cost of goods sold to the operating profit, and decrease the inventory level by the amount of goods sold.
    \item If we were unable to fully meet the demand, we increase the count of stock-out days by one.
    \item We calculate the daily inventory holding costs based on the inventory level at the end of the day, and deduct these costs from the operating profit.
    \item If the current day is the last day of the month, we deduct the fixed monthly costs from the operating profit.
    \item If the current day is not the last day of the simulation, we proceed to step 2.
\end{enumerate}

The simulated operating profits for both SKUs are provided in Table \ref{tab:results_op} (for more detailed data please refer to Appendix A). As observed, for SKU A, the safety stock buffer approach was superior to Models 1 and 2, but Model 3 demonstrated the best performance. This can be explained by the fact that Model 3 properly adjusts order size based on the expected amount of inventory at the time when the order is due. For SKU B, there is one case when the safety stock buffer approach appears to significantly outperform the newsvendor models---the case of uniformly distributed demand. However, since the uniform distribution is clearly inadequate for SKUs with infrequent demand, this outperformance is irrelevant in practice. When considering the triangular or log-normal distributions, we see roughly the same outperformance of Model 3 over the safety stock buffer approach as observed for SKU A, whereas Models 1 and 2 exhibit a substantial drop in profits under the log-normal demand distribution.

\begin{table}[h!]
\caption{Simulated operating profit for different models}\label{tab:results_op}
\renewcommand{\arraystretch}{1.5}
\begin{tabular}{ c R{2.5cm} R{2.5cm} R{2.5cm} R{2.5cm} }
\hline 
\bfseries Distribution &
\bfseries Safety stock &
\bfseries Model 1 &
\bfseries Model 2 &
\bfseries Model 3 \\ 
\hline 
& \multicolumn{4}{c}{SKU A} \\
\hline 
\bfseries Uniform    & 228,553 & 219,307 & 218,205 & 231,235 \\
\bfseries Triangular & 244,149 & 241,121 & 240,062 & 247,676 \\
\bfseries Log-normal & 253,103 & 242,494 & 241,238 & 253,767 \\
\hline 
& \multicolumn{4}{c}{SKU B} \\
\hline 
\bfseries Uniform    &  23,556 &  22,218 &  21,995 &  21,818 \\
\bfseries Triangular &  10,829 &  10,192 &  10,077 &  11,663 \\
\bfseries Log-normal &  10,793 &   8,147 &   7,904 &  10,989 \\
\hline
\end{tabular}
\end{table}

The  most important factor explaining the difference in operating profits is the average inventory on hand (see Table \ref{tab:results_stock}, numbers before the slash). It is evident that inventory levels were significantly higher for the safety stock buffer approach, while all newsvendor models operated with much lower inventory. Consequently, they also experienced a substantial number of stock-out days (see Table \ref{tab:results_stock}, numbers after the slash). Since in our case the inventory holding costs of 2.8 were significantly lower than the gross margin of 40, it was much more beneficial to hold excess inventory rather than risk unmet demand. We can reasonably expect that as holding costs increase, the safety stock buffer approach will begin to lag behind the newsvendor models in performance.

In the previous paragraph we have explained the difference between the safety stock buffer approach and the newsvendor models, but there still remains a significant difference between Models 1-2 and Model 3. The latter experienced fewer stock-outs at about the same inventory levels. This can be attributed to more efficient inventory management in Model 3, which has a much higher base order size (see the $F^{-1}(\cdot)$ arguments in Table \ref{tab:arguments}) but offsets the accumulation of large stock of inventories through dynamic adjustment of the next order size based on the expected inventory level at the date the order is due. In contrast, Models 1 and 2 do not include such an adjustment; they rely on a static order size, which is much smaller, and this is the key reason for their higher stock-out rates. The only notable exception is the uniformly distributed demand for SKU B. However, as previously noted, this case is practically irrevelant, since the uniform distribution is inadequate for SKUs with infrequent demand.

\begin{table}[h!]
\caption{Simulated average inventory and total stock-out amounts for different models}\label{tab:results_stock}
\renewcommand{\arraystretch}{1.5}
\begin{tabular}{ c K{2.5cm} K{2.5cm} K{2.5cm} K{2.5cm} }
\hline 
\bfseries Distribution &
\bfseries Safety stock &
\bfseries Model 1 &
\bfseries Model 2 &
\bfseries Model 3 \\ 
\hline 
& \multicolumn{4}{c}{SKU A} \\
\hline 
\bfseries Uniform    &  4,327 / 0 & 2,052 / 158 & 1,985 / 170 & 2,437 / 34 \\
\bfseries Triangular &  4,194 / 0 & 2,023 / 121 & 1,942 / 136 & 2,136 / 39 \\
\bfseries Log-normal &  4,118 / 0 & 2,007 / 165 & 1,945 / 178 & 2,355 / 49 \\
\hline 
& \multicolumn{4}{c}{SKU B} \\
\hline 
\bfseries Uniform    &  553 / 0 & 198 / 213 & 185 / 236 & 157 / 258 \\
\bfseries Triangular &  661 / 0 & 131 / 261 & 126 / 276 & 174 /  70 \\
\bfseries Log-normal &  661 / 2 & 112 / 406 & 106 / 433 & 225 /  95 \\
\hline
\end{tabular}
\end{table}

\section{Robustness Analysis}

The robustness analysis of Models 1--3 is necessary because we cannot be certain of the true demand distribution, and the distribution assumed when calculating the optimal production batch size in Table \ref{tab:arguments} may differ from the actual demand distribution. The safety stock approach does not rely on any distributional assumptions and should, therefore, be very robust, but it the robustness of the newsvendor models remains an open question that we address here. 

To test the robustness of the results for Models 1--3, we can modify the algorithm described in Section 5 as follows: we continue selecting the optimal production batch size using the original ``model" distribution (i.e., the same values as in Table \ref{tab:arguments}), but at step 5, we simulate demand for the current day with a ``true" distribution, which may or may not be the same as the ``model" distribution.

The results of the robustness analysis, in terms of operating profit, are summarized in Table \ref{tab:robustness_op}. We can immediately observe that the results for the safety stock approach depend solely on the ``true" distribution and are unaffected by the ``model" distribution---that is expected, since the safety stock approach does not utilize any $F^{-1}(\cdot)$ values and thus does not rely on a specific distributional assumption. In contrast, Model 3 results exhibit dynamics quite similar to those of the safety stock approach---the only notable exception is SKU B when the ``true" demand distribution is uniform, but a non-uniform ``model" distribution is assumed. However, as previously noted, this case is unrealistic for SKUs with infrequent demand. Models 1 and 2 appear to be less robust compared to Model 3, especially for SKU B. Nonetheless, it should be noted that the uniform distribution is inadequate for this SKU, so these observations can be considered less relevant.

\begin{table}[h!]
\caption{Operating profits across models and distributions}\label{tab:robustness_op}
\renewcommand{\arraystretch}{1.5}
\begin{tabular}{ K{5cm} R{2cm} R{2cm} R{2cm} R{2cm} }
\hline 
\bfseries Distribution (model/true) &
\bfseries Safety stock &
\bfseries Model 1 &
\bfseries Model 2 &
\bfseries Model 3 \\ 
\hline 
& \multicolumn{4}{c}{SKU A} \\
\hline 
\bfseries Uniform    / Uniform    & 228,553 & 219,307 & 218,205 & 231,235 \\
\bfseries Uniform    / Triangular & 244,060 & 221,089 & 219,818 & 248,234 \\
\bfseries Uniform    / Log-normal & 253,225 & 221,054 & 219,786 & 251,531  \\
\hline
\bfseries Triangular / Uniform    & 228,653 & 205,317 & 209,718 & 230,540 \\
\bfseries Triangular / Triangular & 244,149 & 241,121 & 240,062 & 247,676 \\
\bfseries Triangular / Log-normal & 253,152 & 241,365 & 240,205 & 249,705 \\
\hline
\bfseries Log-normal / Uniform    & 228,519 & 199,791 & 205,023 & 232,001 \\
\bfseries Log-normal / Triangular & 244,066 & 242,053 & 240,971 & 248,606 \\
\bfseries Log-normal / Log-normal & 253,103 & 242,494 & 241,238 & 253,767 \\
\hline 
& \multicolumn{4}{c}{SKU B} \\
\hline 
\bfseries Uniform    / Uniform    &  23,556 &  22,218 &  21,995 &  21,818 \\
\bfseries Uniform    / Triangular &  10,817 & -30,240 & -29,133 &  11,866 \\
\bfseries Uniform    / Log-normal &  10,762 & -30,421 & -29,317 &  11,299 \\
\hline
\bfseries Triangular / Uniform    &  23,575 &  10,444 &  10,315 &  17,036 \\
\bfseries Triangular / Triangular &  10,829 &  10,192 &  10,077 &  11,663 \\
\bfseries Triangular / Log-normal &  10,795 &   9,930 &   9,835 &  10,514 \\
\hline
\bfseries Log-normal / Uniform    &  23,551 &   8,368 &   8,108 &  19,265 \\
\bfseries Log-normal / Triangular &  10,846 &   8,245 &   7,993 &  11,924 \\
\bfseries Log-normal / Log-normal &  10,793 &   8,147 &   7,904 &  10,989 \\
\hline
\end{tabular}
\end{table}

The higher robustness of Model 3 compared to Models 1 and 2, as well as its similarities to the robustness of the safety stock approach, can be explained by the fact that Model 3 adjusts the optimal production batch size based on the expected level of leftover inventory at the time the order is due, whereas Models 1 and 2 do not. The lower robustness of these latter models is related to the mismatch between the ``model" and ``true" demand distributions---particularly between uniform and non-uniform distributions---which leads to either significant excess inventories (resulting in negative profits) or a high number of stock-out days (see Appendix B). While the safety stock buffer approach does not change the order size (which remains constant), it offers flexibility in the timing of placing orders, which is roughly equivalent to having an adjustable batch size with fixed timing.

The robustness analysis also provides three important conclusions for the practical application of the improved newsvendor model when the ``true" demand distribution is unknown:

\begin{enumerate}
    \setlength{\itemsep}{6pt}
    \setlength{\parskip}{0pt}
    \setlength{\parsep}{0pt}
    \item For SKUs with stable demand (such as SKU A in our case), it is generally safer to assume a uniform distribution when determining the optimal production batch size. If using Model 3, a log-normal distribution provides marginally better outcomes, but this is not the case for Models 1 and 2.  
    \item For SKUs with infrequent demand (such as SKU B in our case), Models 1 and 2 should never be used with a uniform distribution, and the triangular distribution appears to be the most robust choice. Model 3 exhibits a quite different behavior: it seems to be robust under the uniform distribution in this case too, despite the theoretical inadequacy of that distribution for such cases; and the second-best choice would be the log-normal distribution.
    \item If only Model 3 is considered and the demand type of the SKU is unknown, using the uniform distribution would be the safest approach.
\end{enumerate}

\section{Conclusions}

In this paper, we propose a multi-period newsvendor model with primary focus on a manufacturing company that needs to determine its production batch size while accounting for inventory leftovers from previous periods. The objective function is to maximize operating profit in the long run, and the model was developed under assumptions consistent with generally accepted accounting principles. In particular, the base case model aligns with the average inventory costing method, and we also discuss potential adjustments for the FIFO method. Inventory holding costs are represented with a single factor; although the model itself allows for differential treatment of warehousing costs, shrinkage costs, and the cost of working capital, all of these could generally be incorporated into a single variable.

We used Monte Carlo simulation to compare the outcomes of the proposed model with those of the original newsvendor model, a simple multi-period extension, and a safety stock buffer approach. The simulation was conducted on two SKUs of a real-world clothing manufacturing: one SKU with relatively stable daily demand and another with unstable, infrequent demand. We used three ``model" distributions---uniform, triangular, and log-normal---for calculating the optimal batch size and modeling daily demand. The simulation demonstrated that the proposed model generally results in higher operating profits than other newsvendor models and, in most cases, outperforms the safety stock buffer approach. We also examined the factors influencing these outcomes, and found that the proposed model results in fewer stock-outs compared to other newsvendor models and maintains lower inventory levels than the safety stock buffer approach. While the proposed model is not free of stock-outs, and the lost profits due to unmet demand remain an issue, it demonstrated better adaptability in adjusting production batch size relative to the other models analyzed.

Additionally, we conducted a robustness analysis by introducing a mismatch between the ``model" and the ``true" demand distribution. As expected, the safety stock approach proved to be the most robust, since it is distribution-free. Classic newsvendor models showed limited robustness, failing when the ``true" demand distribution was uniform and the ``model" one was triangular or log-normal, which resulted in dramatic operating losses in the case of unstable, infrequent demand. However, the proposed model, due to its better adaptability in adjusting the production batch size, proved to be nearly as robust as the safety stock buffer approach. Although some underperformance was observed when the ``true" distribution was uniform and demand was unstable and infrequent, the model still generated solid operating profits, and in most cases, the mismatch between the ``model" and the ``true" distribution did not cause critical drops in profitability. Finally, we provided recommendation regarding the choice of the ``model" distribution for determining the optimal production batch size, depending on whether the SKUs' demand patterns are known or unknown.

\section*{Declarations}

\paragraph{Funding.} The author did not receive support from any organization for the submitted work.

\paragraph{Competing interests.} The author has no financial or proprietary interests in any material discussed in this article.

\paragraph{AI Usage.} DeepAI (deepai.org) with Standard model was used solely for grammar correction and language editing. All scientific content, interpretations, and conclusions were generated by the author. The final text was reviewed and verified by the author.

\bibliographystyle{plainnat}
 \bibliography{main}

\newpage

\begin{appendices}

\appendix

\section{Detailed Simulation Results}
\setcounter{table}{0}
\renewcommand{\thetable}{A\arabic{table}}

Note: In the tables below \textit{MOE 95\%} specifies the margin of error at the 95\% confidence level for the metric above it; \textbf{\textit{x}th perc.} refers to the \textbf{\textit{x}}-th percentile.

\renewcommand{\arraystretch}{1.145}
\begin{longtable}{ L{3cm} L{2cm} R{1.8cm} R{1.8cm} R{1.8cm} R{1.8cm} }\caption{Detailed simulation results for SKU A}\label{tab:detailed_A}\\
\toprule
\textbf{Performance indicator} & \textbf{Metric} & \textbf{Safety stock} & \textbf{Model 1} & \textbf{Model 2} & \textbf{Model 3} \\ 
\midrule
\endfirsthead

\multicolumn{2}{c}{{Continued from previous page...}{}} \\
\toprule
\textbf{Performance indicator} & \textbf{Metric} & \textbf{Safety stock} & \textbf{Model 1} & \textbf{Model 2} & \textbf{Model 3} \\ 
\midrule
\endhead

\midrule
\multicolumn{6}{r}{{Continued on next page...}} \\
\endfoot

\bottomrule
\endlastfoot
& \multicolumn{5}{c}{SKU A, ``model" and ``true" distribution: Uniform} \\
\hline 
\bfseries \multirow{6}{2.5cm}{\centering Operating profit}
 &\bfseries Average    & 228,553 & 219,307 & 218,205 & 231,235 \\
 &\textit{MOE 95\%}     &   ±198  &     ±46 &     ±40 &    ±179 \\
 &\bfseries St.dev.    &  3032.6 &   709.0 &   615.4 &  2733.7 \\
 &\bfseries Median     & 228,499 & 219,469 & 218,339 & 231,182 \\
 &\bfseries 10th perc. & 224,713 & 218,436 & 217,446 & 227,824 \\
 &\bfseries  5th perc. & 223,494 & 218,063 & 217,089 & 226,738 \\
\hline 
\bfseries \multirow{6}{2.5cm}{\centering Average inventory on hand}
 &\bfseries Average    &   4,327 &   2,052 &   1,985 &   2,427 \\
 &\textit{MOE 95\%}     &    ±1.9 &   ±14.5 &   ±12.5 &    ±2.2 \\
 &\bfseries St.dev.    &    29.4 &   221.8 &   190.9 &    33.1 \\
 &\bfseries Median     &   4,327 &   2,002 &   1,943 &   2,427 \\
 &\bfseries 95th perc. &   4,375 &   2,467 &   2,339 &   2,483 \\
 &\bfseries 99th perc. &   4,394 &   2,883 &   2,679 &   2,503 \\
\hline
\bfseries \multirow{6}{2.5cm}{\centering Stock-out days}
 &\bfseries Average    &       0	&     158 &     170 &      34 \\
 &\textit{MOE 95\%}     &    ±0.0 &    ±1.8 &    ±1.8 &    ±0.4 \\
 &\bfseries St.dev.    &     0.1 &    27.3 &    27.2 &     6.5 \\
 &\bfseries Median     &       0	&     158 &     169 &      34 \\
 &\bfseries 95th perc. &       0	&     204 &     214 &      44 \\
 &\bfseries 99th perc. &       0	&     220 &     233 &      51 \\
\hline
& \multicolumn{5}{c}{SKU A, ``model" and ``true" distribution: Triangular} \\
\hline 
\bfseries \multirow{6}{2.5cm}{\centering Operating profit}
 &\bfseries Average    & 244,149 & 241,121 & 240,062 & 247,676 \\
 &\textit{MOE 95\%}     &   ±132  &     ±38 &     ±30 &    ±116 \\
 &\bfseries St.dev.    &  2016.0 &   584.8 &   460.2 &  1776.7 \\
 &\bfseries Median     & 244,057 & 241,249 & 240,159 & 247,660 \\
 &\bfseries 10th perc. & 241,594 & 240,412 & 239,506 & 245,370 \\
 &\bfseries  5th perc. & 240,885 & 240,047 & 239,234 & 244,719 \\
\pagebreak
\bfseries \multirow{6}{2.5cm}{\centering Average inventory on hand}
 &\bfseries Average    &   4,194 &   2,023 &   1,942 &   2,136 \\
 &\textit{MOE 95\%}     &    ±1.4 &   ±11.7 &    ±9.1 &    ±1.4 \\
 &\bfseries St.dev.    &    20.8 &   179.8 &   138.9 &    21.9 \\
 &\bfseries Median     &   4,195 &   1,985 &   1,914 &   2,136 \\
 &\bfseries 95th perc. &   4,226 &   2,364 &   2,219 &   2,172 \\
 &\bfseries 99th perc. &   4,242 &   2,632 &   2,380 &   2,186 \\
\hline
\bfseries \multirow{6}{2.5cm}{\centering Stock-out days}
 &\bfseries Average    &       0	&     121 &     136 &      39 \\
 &\textit{MOE 95\%}     &    ±0.0 &    ±1.6 &    ±1.5 &    ±0.4 \\
 &\bfseries St.dev.    &     0.0 &    24.0 &    23.2 &     6.5 \\
 &\bfseries Median     &       0	&     122 &     137 &      39 \\
 &\bfseries 95th perc. &       0	&     161 &     174 &      51 \\
 &\bfseries 99th perc. &       0	&     177 &     189 &      55 \\
\hline
& \multicolumn{5}{c}{SKU A, ``model" and ``true" distribution: Log-normal} \\
\hline 
\bfseries \multirow{6}{2.5cm}{\centering Operating profit}
 &\bfseries Average    & 253,103 & 242,494 & 241,238 & 253,767 \\
 &\textit{MOE 95\%}     &   ±181  &     ±34 &     ±29 &    ±153 \\
 &\bfseries St.dev.    &  2775.8	&   515.1 &   450.1 &  2341.4 \\
 &\bfseries Median     & 253,099 & 242,583 & 241,310 & 253,811 \\
 &\bfseries 10th perc. & 249,635 & 241,865 & 240,688 & 250,712 \\
 &\bfseries  5th perc. & 248,643 & 241,579 & 240,462 & 249,837 \\
\hline 
\bfseries \multirow{6}{2.5cm}{\centering Average inventory on hand}
 &\bfseries Average    &   4,118 &   2,007 &   1,945 &   2,355 \\
 &\textit{MOE 95\%}     &    ±1.8 &    ±9.8 &    ±8.4 &    ±1.9 \\
 &\bfseries St.dev.    &    27.2	&   149.5 &   128.4 &    29.2 \\
 &\bfseries Median     &   4,118 &   1,986 &   1,926 &   2,356 \\
 &\bfseries 95th perc. &   4,163 &   2,272 &   2,168 &   2,402 \\
 &\bfseries 99th perc. &   4,181 &   2,536 &   2,399 &   2,420 \\
\hline
\bfseries \multirow{6}{2.5cm}{\centering Stock-out days}
 &\bfseries Average    &       0	&     165 &     178 &      49 \\
 &\textit{MOE 95\%}     &    ±0.0 &    ±1.6 &    ±1.6 &    ±0.5 \\
 &\bfseries St.dev.    &     0.7 &    24.9 &    24.7 &     7.8 \\
 &\bfseries Median     &       0	&     164 &     177 &      49 \\
 &\bfseries 95th perc. &       2	&     204 &     217 &      63 \\
 &\bfseries 99th perc. &       3	&     219 &     234 &      69 \\
\hline
\end{longtable}

\renewcommand{\arraystretch}{1.145}
\begin{longtable}{ L{3cm} L{2cm} R{1.8cm} R{1.8cm} R{1.8cm} R{1.8cm} }\caption{Detailed simulation results for SKU B}\label{tab:detailed_A}\\
\toprule
\textbf{Performance indicator} & \textbf{Metric} & \textbf{Safety stock} & \textbf{Model 1} & \textbf{Model 2} & \textbf{Model 3} \\ 
\midrule
\endfirsthead

\multicolumn{2}{c}{{Continued from previous page...}} \\
\toprule
\textbf{Performance indicator} & \textbf{Metric} & \textbf{Safety stock} & \textbf{Model 1} & \textbf{Model 2} & \textbf{Model 3} \\ 
\midrule
\endhead

\midrule
\multicolumn{6}{r}{{Continued on next page...}} \\
\endfoot

\bottomrule
\endlastfoot
& \multicolumn{5}{c}{SKU B, ``model" and ``true" distribution: Uniform} \\
\hline 
\bfseries \multirow{6}{2.5cm}{\centering Operating profit}
 &\bfseries Average    &  23,556 &  22,218 &  21,995 &  21,818 \\
 &\textit{MOE 95\%}     &    ±28  &      ±7 &      ±5 &     ±16 \\
 &\bfseries St.dev.    &   430.6 &   101.0 &    83.8 &   242.8 \\
 &\bfseries Median     &  23,545 &  22,235 &  22,010 &  21,815 \\
 &\bfseries 10th perc. &  23,036 &  22,093 &  21,887 &  21,512 \\
 &\bfseries  5th perc. &  22,844 &  22,032 &  21,835 &  21,419 \\
\hline 
\bfseries \multirow{6}{2.5cm}{\centering Average inventory on hand}
 &\bfseries Average    &     553 &     198 &     185 &     157 \\
 &\textit{MOE 95\%}     &    ±0.3 &    ±2.0 &    ±1.7 &    ±0.2 \\
 &\bfseries St.dev.    &     4.2 &    31.0 &    25.4 &     3.3 \\
 &\bfseries Median     &     553 &     192 &     180 &     157 \\
 &\bfseries 95th perc. &     560 &     256 &     231 &     163 \\
 &\bfseries 99th perc. &     563 &     297 &     267 &     165 \\
\hline
\bfseries \multirow{6}{2.5cm}{\centering Stock-out days}
 &\bfseries Average    &       0	&     213 &     236 &     256 \\
 &\textit{MOE 95\%}     &    ±0.0 &    ±2.4 &    ±2.3 &    ±1.3 \\
 &\bfseries St.dev.    &     0.0 &    36.2 &    35.8 &    19.5 \\
 &\bfseries Median     &       0	&     213 &     236 &     259 \\
 &\bfseries 95th perc. &       0	&     271 &     292 &     289 \\
 &\bfseries 99th perc. &       0	&     292 &     314 &     301 \\
\hline
& \multicolumn{5}{c}{SKU B, ``model" and ``true" distribution: Triangular} \\
\hline 
\bfseries \multirow{6}{2.5cm}{\centering Operating profit}
 &\bfseries Average    &  10,829 &  10,192 &  10,077 &  11,663 \\
 &\textit{MOE 95\%}     &    ±22  &      ±4 &      ±4 &     ±19 \\
 &\bfseries St.dev.    &   340.2 &    59.2 &    53.6 &   286.9 \\
 &\bfseries Median     &  10,829 &  10,200 &  10,084 &  11,663 \\
 &\bfseries 10th perc. &  10,408 &  10,120 &  10,011 &  11,302 \\
 &\bfseries  5th perc. &  10,269 &  10,093 &   9,988 &  11,191 \\
\pagebreak
\bfseries \multirow{6}{2.5cm}{\centering Average inventory on hand}
 &\bfseries Average    &     661 &     131 &     126 &     174 \\
 &\textit{MOE 95\%}     &    ±0.2 &    ±1.1 &    ±1.0 &    ±0.2 \\
 &\bfseries St.dev.    &     3.6	&    17.3 &    15.5 &     3.5 \\
 &\bfseries Median     &     661 &     128 &     123 &     174 \\
 &\bfseries 95th perc. &     667 &     163 &     154 &     180 \\
 &\bfseries 99th perc. &     669 &     185 &     175 &     182 \\
\hline
\bfseries \multirow{6}{2.5cm}{\centering Stock-out days}
 &\bfseries Average    &       0	&     261 &     276 &      70 \\
 &\textit{MOE 95\%}     &    ±0.0 &    ±2.5 &    ±2.5 &    ±0.7 \\
 &\bfseries St.dev.    &     0.0 &    38.1 &    37.8 &    11.2 \\
 &\bfseries Median     &       0	&     260 &     276 &      69 \\
 &\bfseries 95th perc. &       0	&     324 &     339 &      88 \\
 &\bfseries 99th perc. &       0	&     350 &     364 &      96 \\
\hline
& \multicolumn{5}{c}{SKU B, ``model" and ``true" distribution: Log-normal} \\
\hline 
\bfseries \multirow{6}{2.5cm}{\centering Operating profit}
 &\bfseries Average    &  10,793 &   8,147 &   7,904 &  10,989 \\
 &\textit{MOE 95\%}     &    ±37  &      ±3 &      ±3 &     ±27 \\
 &\bfseries St.dev.    &   560.4 &	 45.5 &    40.4 &   419.9 \\
 &\bfseries Median     &  10,777 &   8,154 &   7,911 &  11,001 \\
 &\bfseries 10th perc. &  10,079 &   8,091 &   7,855 &  10,458 \\
 &\bfseries  5th perc. &   9,890 &   8,062 &   7,828 &  10,303 \\
\hline 
\bfseries \multirow{6}{2.5cm}{\centering Average inventory on hand}
 &\bfseries Average    &     661 &     112 &     106 &     225 \\
 &\textit{MOE 95\%}     &    ±0.4 &    ±0.8 &    ±0.7 &    ±0.3 \\
 &\bfseries St.dev.    &     5.6 &    12.1 &    10.5 &     5.2 \\
 &\bfseries Median     &     661 &     111 &     104 &     225 \\
 &\bfseries 95th perc. &     670 &     134 &     125 &     233 \\
 &\bfseries 99th perc. &     673 &     149 &     137 &     236 \\
\hline
\bfseries \multirow{6}{2.5cm}{\centering Stock-out days}
 &\bfseries Average    &       2	&     406 &     433 &      95 \\
 &\textit{MOE 95\%}     &    ±0.2 &    ±3.1 &    ±3.1 &    ±1.2 \\
 &\bfseries St.dev.    &     2.8 &    46.7 &    46.8 &    18.9 \\
 &\bfseries Median     &       0	&     406 &     433 &      94 \\
 &\bfseries 95th perc. &       8	&     481 &     509 &     128 \\
 &\bfseries 99th perc. &      12	&     506 &     534 &     145 \\
\hline
\end{longtable}

\section{Robustness Analysis Results}
\setcounter{table}{0}
\renewcommand{\thetable}{B\arabic{table}}

\renewcommand{\arraystretch}{1.19}

\begin{longtable}{ L{3cm} L{2cm} R{1.8cm} R{1.8cm} R{1.8cm} R{1.8cm} }\caption{Robustness analysis results for SKU A}\label{tab:detailed_A}\\
\toprule
\textbf{Distribution, model / true} & \textbf{Metric} & \textbf{Safety stock} & \textbf{Model 1} & \textbf{Model 2} & \textbf{Model 3} \\ 
\midrule
\endfirsthead

\multicolumn{2}{c}{{Continued from previous page...}} \\
\toprule
\textbf{Distribution, model / true} & \textbf{Metric} & \textbf{Safety stock} & \textbf{Model 1} & \textbf{Model 2} & \textbf{Model 3} \\ 
\midrule
\endhead

\midrule
\multicolumn{6}{r}{{Continued on next page...}} \\
\endfoot

\bottomrule
\endlastfoot
& \multicolumn{5}{c}{SKU A, Operating profit} \\
\hline 
\bfseries \multirow{2}{2.5cm}{\centering Uniform / Triangular}
 &\bfseries Average    & 244,060 & 221,089 & 219,818 & 248,234 \\
 &\bfseries St.dev.    &  1995.4 &    96.0 &    91.6 &  1862.9 \\
\hline
\bfseries \multirow{2}{2.5cm}{\centering Uniform / Log-normal}
 &\bfseries Average    & 253,225 & 221,054 & 219,786 & 251,531 \\
 &\bfseries St.dev.    &  2961.7 &   136.4 &   131.5 &  2277.6 \\
\hline
\bfseries \multirow{2}{2.5cm}{\centering Triangular / Uniform}
 &\bfseries Average    & 228,653 & 205,317 & 209,718 & 230,540 \\
 &\bfseries St.dev.    &  2982.5 & 14163.0 & 13387.4 &  2578.7 \\
\hline
\bfseries \multirow{2}{2.5cm}{\centering Triangular / Log-normal}
 &\bfseries Average    & 253,152 & 241,365 & 240,205 & 249,705 \\
 &\bfseries St.dev.    &  2883.9 &   460.8 &   404.8 &  2019.1 \\
\hline 
\bfseries \multirow{2}{2.5cm}{\centering Log-normal / Uniform}
 &\bfseries Average    & 228,519 & 199,791 & 205,023 & 232,001 \\
 &\bfseries St.dev.    &  3040.3 & 14786.6 & 14220.3 &  2866.2 \\
\hline
\bfseries \multirow{2}{2.5cm}{\centering Log-normal / Triangular}
 &\bfseries Average    & 244,066 & 242,053 & 240,971 & 248,606 \\
 &\bfseries St.dev.    &  1943.8 &   761.5 &   584.3 &  1903.4 \\
\hline
& \multicolumn{5}{c}{SKU A, Average inventory on hand} \\
\hline 
\bfseries \multirow{2}{2.5cm}{\centering Uniform / Triangular}
 &\bfseries Average    &   4,195 &   1,488 &   1,474 &   2,234 \\
 &\bfseries St.dev.    &    21.1 &    21.3 &    19.8 &    22.3 \\
\hline
\bfseries \multirow{2}{2.5cm}{\centering Uniform / Log-normal}
 &\bfseries Average    &   4,118 &   1,498 &   1,483 &   2,178 \\
 &\bfseries St.dev.    &    28.4 &    29.8 &    28.0 &    28.8 \\
\hline
\bfseries \multirow{2}{2.5cm}{\centering Triangular / Uniform}
 &\bfseries Average    &   4,325 &  12,449 &  10,819 &   2,330 \\
 &\bfseries St.dev.    &    28.3 &  4217.0 &  3980.0 &    31.4 \\
\hline
\bfseries \multirow{2}{2.5cm}{\centering Triangular / Log-normal}
 &\bfseries Average    &   4,118 &   1,947 &   1,897 &   2,092 \\
 &\bfseries St.dev.    &    28.1	&   136.0 &   117.0 &    26.3 \\
\hline 
\bfseries \multirow{2}{2.5cm}{\centering Log-normal / Uniform}
 &\bfseries Average    &   4,327 &  14,417 &  12,508 &   2,615 \\
 &\bfseries St.dev.    &    29.4 &  4390.5 &  4214.7 &    34.8 \\
\hline
\bfseries \multirow{2}{2.5cm}{\centering Log-normal / Triangular}
 &\bfseries Average    &   4,194 &   2,143 &   2,028 &   2,426 \\
 &\bfseries St.dev.    &    19.6 &   235.0 &   178.7 &    22.7 \\

& \multicolumn{5}{c}{SKU A, Stock-out days} \\
\hline 
\bfseries \multirow{2}{2.5cm}{\centering Uniform / Triangular}
 &\bfseries Average    &       0 &     320 &     330 &      25 \\
 &\bfseries St.dev.    &     0.0 &    16.6 &    16.2 &     5.2 \\
\hline
\bfseries \multirow{2}{2.5cm}{\centering Uniform / Log-normal}
 &\bfseries Average    &       1 &     351 &     361 &      78 \\
 &\bfseries St.dev.    &     0.8 &    21.4 &    21.3 &    10.2 \\
\hline
\bfseries \multirow{2}{2.5cm}{\centering Triangular / Uniform}
 &\bfseries Average    &       0 &       7 &       9 &      45 \\
 &\bfseries St.dev.    &     0.1 &     7.3 &     8.8 &     7.5 \\
\hline
\bfseries \multirow{2}{2.5cm}{\centering Triangular / Log-normal}
 &\bfseries Average    &       1 &     177 &     189 &      97 \\
 &\bfseries St.dev.    &     0.8 &    24.8 &    24.4 &    11.2 \\
\hline 
\bfseries \multirow{2}{2.5cm}{\centering Log-normal / Uniform}
 &\bfseries Average    &       0 &       5 &       7 &      18 \\
 &\bfseries St.dev.    &     0.1 &     6.2 &     7.6 &     4.8 \\
\hline
\bfseries \multirow{2}{2.5cm}{\centering Log-normal / Triangular}
 &\bfseries Average    &       0 &     105 &     121 &       9 \\
 &\bfseries St.dev.    &     0.0 &    23.1 &    22.5 &     3.1 \\
\hline
\end{longtable}

\begin{longtable}{ L{3cm} L{2cm} R{1.8cm} R{1.8cm} R{1.8cm} R{1.8cm} }\caption{Robustness analysis results for SKU B}\label{tab:detailed_A}\\
\toprule
\textbf{Distribution, model / true} & \textbf{Metric} & \textbf{Safety stock} & \textbf{Model 1} & \textbf{Model 2} & \textbf{Model 3} \\ 
\midrule
\endfirsthead

\multicolumn{2}{c}{{Continued from previous page...}} \\
\toprule
\textbf{Distribution, model / true} & \textbf{Metric} & \textbf{Safety stock} & \textbf{Model 1} & \textbf{Model 2} & \textbf{Model 3} \\ 
\midrule
\endhead

\midrule
\multicolumn{6}{r}{{Continued on next page...}} \\
\endfoot

\bottomrule
\endlastfoot
& \multicolumn{5}{c}{SKU B, Operating profit} \\
\hline 
\bfseries \multirow{2}{2.5cm}{\centering Uniform / Triangular}
 &\bfseries Average    &  10,817 & -30,240 & -29,133 &  11,866 \\
 &\bfseries St.dev.    &   356.0 &  1964.7 &  1964.5 &   355.5 \\
\hline
\bfseries \multirow{2}{2.5cm}{\centering Uniform / Log-normal}
 &\bfseries Average    &  10,762 & -30,421 & -29,317 &  11,299 \\
 &\bfseries St.dev.    &   540.7 &  3169.2 &  3168.0 &   445.4 \\
\hline
\bfseries \multirow{2}{2.5cm}{\centering Triangular / Uniform}
 &\bfseries Average    &  23,575 &  10,444 &  10,315 &  17,036 \\
 &\bfseries St.dev.    &   423.7 &    12.4 &    12.3 &   110.2 \\
\hline
\bfseries \multirow{2}{2.5cm}{\centering Triangular / Log-normal}
 &\bfseries Average    &  10,795 &   9,930 &   9,835 &  10,514 \\
 &\bfseries St.dev.    &   551.5 &   168.5 &   149.7 &   354.7 \\
\hline 
\bfseries \multirow{2}{2.5cm}{\centering Log-normal / Uniform}
 &\bfseries Average    &  23,551 &   8,368 &   8,108 &  19,265 \\
 &\bfseries St.dev.    &   438.6 &    11.9 &    11.9 &   159.7 \\

\bfseries \multirow{2}{2.5cm}{\centering Log-normal / Triangular}
 &\bfseries Average    &  10,846 &   8,245 &   7,993 &  11,924 \\
 &\bfseries St.dev.    &   354.1 &    19.7 &    18.1 &   331.3 \\
\hline
& \multicolumn{5}{c}{SKU B, Average inventory on hand} \\
\hline 
\bfseries \multirow{2}{2.5cm}{\centering Uniform / Triangular}
 &\bfseries Average    &     661 &  15,324 &  14,929 &     279 \\
 &\bfseries St.dev.    &     3.9 &   592.2 &   592.2 &     4.2 \\
\hline
\bfseries \multirow{2}{2.5cm}{\centering Uniform / Log-normal}
 &\bfseries Average    &     662 &  15,377 &  14,983 &     285 \\
 &\bfseries St.dev.    &     5.5 &   961.3 &   960.9 &     5.3 \\
\hline
\bfseries \multirow{2}{2.5cm}{\centering Triangular / Uniform}
 &\bfseries Average    &     553 &      51 &      51 &      95 \\
 &\bfseries St.dev.    &     4.1 &     1.2 &     1.2 &     2.0 \\
\hline
\bfseries \multirow{2}{2.5cm}{\centering Triangular / Log-normal}
 &\bfseries Average    &     661 &     216 &     204 &     186 \\
 &\bfseries St.dev.    &     5.7 &    53.5 &    47.3 &     4.5 \\
\hline 
\bfseries \multirow{2}{2.5cm}{\centering Log-normal / Uniform}
 &\bfseries Average    &     553 &      42 &      41 &     117 \\
 &\bfseries St.dev.    &     4.3 &     1.0 &     1.0 &     2.5 \\
\hline
\bfseries \multirow{2}{2.5cm}{\centering Log-normal / Triangular}
 &\bfseries Average    &     661 &      81 &      77 &     215 \\
 &\bfseries St.dev.    &     3.9 &     4.8 &     4.3 &     4.0 \\
\hline
& \multicolumn{5}{c}{SKU B, Stock-out days} \\
\hline 
\bfseries \multirow{2}{2.5cm}{\centering Uniform / Triangular}
 &\bfseries Average    &       0 &       0 &       0 &       3 \\
 &\bfseries St.dev.    &     0.0 &     0.4 &     0.4 &     2.2 \\
\hline
\bfseries \multirow{2}{2.5cm}{\centering Uniform / Log-normal}
 &\bfseries Average    &       2 &       1 &       1 &      47 \\
 &\bfseries St.dev.    &     2.9 &     1.5 &     1.5 &    13.6 \\
\hline
\bfseries \multirow{2}{2.5cm}{\centering Triangular / Uniform}
 &\bfseries Average    &       0 &   1,145 &   1,155 &     651 \\
 &\bfseries St.dev.    &     0.0 &    25.6 &    25.5 &    24.4 \\
\hline
\bfseries \multirow{2}{2.5cm}{\centering Triangular / Log-normal}
 &\bfseries Average    &       2 &     201 &     214 &     151 \\
 &\bfseries St.dev.    &     2.7 &    47.1 &    46.9 &    22.5 \\
\hline 
\bfseries \multirow{2}{2.5cm}{\centering Log-normal / Uniform}
 &\bfseries Average    &       0 &   1,291 &   1,309 &     472 \\
 &\bfseries St.dev.    &     0.0 &    24.6 &    24.5 &    23.9 \\
\hline
\bfseries \multirow{2}{2.5cm}{\centering Log-normal / Triangular}
 &\bfseries Average    &       0 &     503 &     532 &      23 \\
 &\bfseries St.dev.    &     0.0 &    37.2 &    36.8 &     6.7 \\
\hline
\end{longtable}

\end{appendices}

\end{document}